\documentclass[11pt, reqno]{amsart}

\usepackage{amsmath,amssymb,amsthm, epsfig}
\usepackage{mathalfa}
\usepackage{amsfonts}
\usepackage[mathscr]{euscript}
\usepackage{wrapfig}
\usepackage{graphicx}
\usepackage{color}

\graphicspath{ {images/} }

\def \N {\mathbb{N}}
\def \R {\mathbb{R}}

\numberwithin{equation}{section}

\newcommand{\intav}[1]{\mathchoice {\mathop{\vrule width 6pt height 3 pt depth  -2.5pt
\kern -8pt \intop}\nolimits_{\kern -6pt#1}} {\mathop{\vrule width
5pt height 3  pt depth -2.6pt \kern -6pt \intop}\nolimits_{#1}}
{\mathop{\vrule width 5pt height 3 pt depth -2.6pt \kern -6pt
\intop}\nolimits_{#1}} {\mathop{\vrule width 5pt height 3 pt depth
-2.6pt \kern -6pt \intop}\nolimits_{#1}}}

\title[Sharp regularity for degenerate pdes]{Geometric tangential analysis and sharp regularity for degenerate pdes}

\author[E.V.~Teixeira]{Eduardo V. Teixeira}
\address{University of Central Florida, 4393 Andromeda Loop N, Orlando, FL 32816, USA}{}
\email{eduardo.teixeira@ucf.edu}

\author[J.M.~Urbano]{Jos\'{e} Miguel Urbano}
\address{CMUC, Department of Mathematics, University of Coimbra, 3001-501 Coimbra, Portugal}{}
\email{jmurb@mat.uc.pt}

\begin{document}

\begin{abstract} 
We provide a broad overview on qualitative versus quantitative regularity estimates in the theory of degenerate parabolic pdes. The former relates to DiBenedetto's revolutionary method of intrinsic scaling, while the latter is achieved by means of what has been termed geometric tangential analysis. We discuss, in particular, sharp estimates for the parabolic $p-$Poisson equation, for the porous medium equation and for the doubly nonlinear equation.
  
\end{abstract}

\dedicatory{To Emmanuele DiBenedetto, with admiration and friendship,\\ on the occasion of his 70th birthday.}

\date{\today}

\maketitle

\section{Introduction}

In the mid 1980's, Emmanuele DiBenedetto made a series of crucial contributions (see, e.g., \cite{DiBe82, DiBe831, DiBe832, DiBeFr84, DiBeFr85}) to the understanding of the regularity properties of weak solutions of singular and degenerate parabolic equations. His \textit{method of intrinsic scaling} (cf. \cite{DiBe93, DiBeGiaVes11, DiBeUrbVes04, Urb08} for rather complete accounts) would become a landmark in regularity theory, to be used extensively in the next decades to treat a variety of pdes, the most celebrated being the $p-$Laplace equation and the porous medium equation. The main insight supporting the method is that each degenerate pde must be analysed in its own geometric setting, in which space and time scale according to the nature of the degeneracy. The concrete implementation of this general principle is rather involved, requiring the use of fine analytic estimates and sophisticated iterative methods (see, among others, \cite{DiBeChen, DiBeGiaVes081, DiBeGiaVes082}). Although powerful and versatile, the method of intrinsic scaling delivers essentially qualitative results, placing the weak solutions of certain pdes in the right regularity class but not providing sharp quantitative information. For example, it tells us that solutions of the porous medium equation 
$$u_{t} - {\rm div} \left( m |u|^{m-1} \nabla u \right) = f$$
are locally of class $C^{0, \alpha}$, for a certain (small) H\"older exponent $\alpha$ but is tight-lipped on the best possible $\alpha$.

While, for many purposes, qualitative estimates for a given model are enough, the investigation pertaining to sharp estimates in diffusive pdes is, by no means, a mere fanciful inquire. On the contrary, sharp estimates reveal important nuances of the problem and play a decisive role in a finer analysis of the model. As a way of example, they  are decisive in the investigation of problems involving free boundaries. As a structural attribute of a given pde, obtaining optimal regularity estimates for a given diffusive model is often a challenging problem. In what follows, we will describe a successful geometric approach which often leads to sharp, or at least improved, estimates in H\"older spaces for solutions of non-homogeneous pdes; the method is inspired by geometric insights related to the notion of \textit{tangent pde models}. We emphasise that our approach and the method of intrinsic scaling are not in competition but rather complement each other. 

In the sequel, we will describe in more detail a few seminal ideas which foster intuitive insights leading to a technical apparatus supporting the method. We will also  exemplify the power of those ideas in some concrete, relevant problems, namely in obtaining sharp estimates for the parabolic $p-$Poisson equation, for the porous medium equation and for pdes involving doubly nonlinearities.

\section{Geometric tangential analysis}

The abstract  concept of Tangent is rather classical and widely spread in the realm of mathematical sciences. It bears a notion of approximation, usually involving more regular objects, from which one can infer pertinent information about the original entity. Probably one of the most well known examples of Tangent comes from the idea of differentiation, where one locally approximates a nonlinear map by a linear one. The acclaimed Inverse Function Theorem from Calculus asserts that if the linear approximation of a function $f$ at a point $a$ is injective and surjective, then so is $f$ in a neighbourhood of  $a$. This is a classical example where qualitative information on the approximating object is transferred to the approximated one.  

While tangent lines to the graph of a function, or even tangent hyperplanes to manifolds are, \textit{per se}, rather concrete manifestations of Tangent, this powerful mathematical concept transcends to more abstract settings, ultimately yielding decisive breakthroughs. 
 
In the lines of the analogy above, Geometric Tangential Analysis (GTA) refers to a constructive systematic approach based on the concept that a problem which enjoys greater regularity can be tangentially accessed by certain classes of pdes. By means of iterative arguments, the method then imports this regularity, properly corrected through the path used to access the tangential equation, to the original class. The roots of this idea likely go back to the foundation of De Giorgi's geometric measure theory of minimal surfaces, and accordingly, it is present in the development of the modern theory of free boundary problems. 

Indeed, an instrumental argument in De Giorgi's geometric measure treatment of minimal surfaces is the so called \textit{flatness improvement}. Roughly speaking, De Giorgi in \cite{DGi} shows that if a minimal surface $S$ is flat enough, say in $B_1$, then it is even flatter in $B_{1/2}$. At least as important as the theorem itself is the reasoning of its proof, which heuristically goes as follows: arguing by contradiction, one would produce a sequence of minimal surfaces $S_j$ in $B_1$, that are $1/j$ flat with respect to a direction $\nu_j$; however the aimed flatness improvement is not verified in $B_{1/2}$. By {compactness} arguments, an appropriate scaling of $S_j$ converges to the graph of a function $f$, which ought to solve the linearised equation, namely $\Delta f  =0$.  Since the limiting function $f$ is very smooth, flatness improvement is verified for $f$. Thus, for $j_0$ sufficiently large, one reaches a contraction on the assumption that no flatness improvement was possible for $S_{j_0}$. 

Such a revolutionary, seminal idea borne fruit in many other fields of research. In particular the motto \textit{flatness implies regularity}, largely promoted by Caffarelli and collaborators,  thrived in the theory of free boundary problems from the 1970's and 1980's (see, among others, \cite{AC, ACF, Caff77, Caff80, CC1, CC2, CC3}). Powerful methods and geometric insights designed for the study of free boundary problems evolved and, in the 1990's, played a decisive role in Caffarelli's work on fully non-linear elliptic pdes (cf. \cite{C1}) and, subsequently, in his studies on Monge-Amp\`ere equations (see \cite{CMA}). As for second order fully nonlinear elliptic equations, Caffarelli uses Krylov-Safonov Harnack inequality, designed for viscosity solutions, as a universal compactness device. He measures \textit{closeness} between variable coefficient equations and constant coefficient equations by means of coefficient oscillation; no linearisation takes place. Ultimately, he shows that if the constant coefficient equation $F(x_0, D^2u)=0$ has a good regularity theory, then $F(x, D^2u)=0$ inherits some universal estimates, provided the  coefficient oscillation is small enough. 

Restricted to diffusive processes, perhaps a didactical way to contemplate GTA is by drawing connected dots. 
\begin{center}
    \includegraphics[scale=0.07]{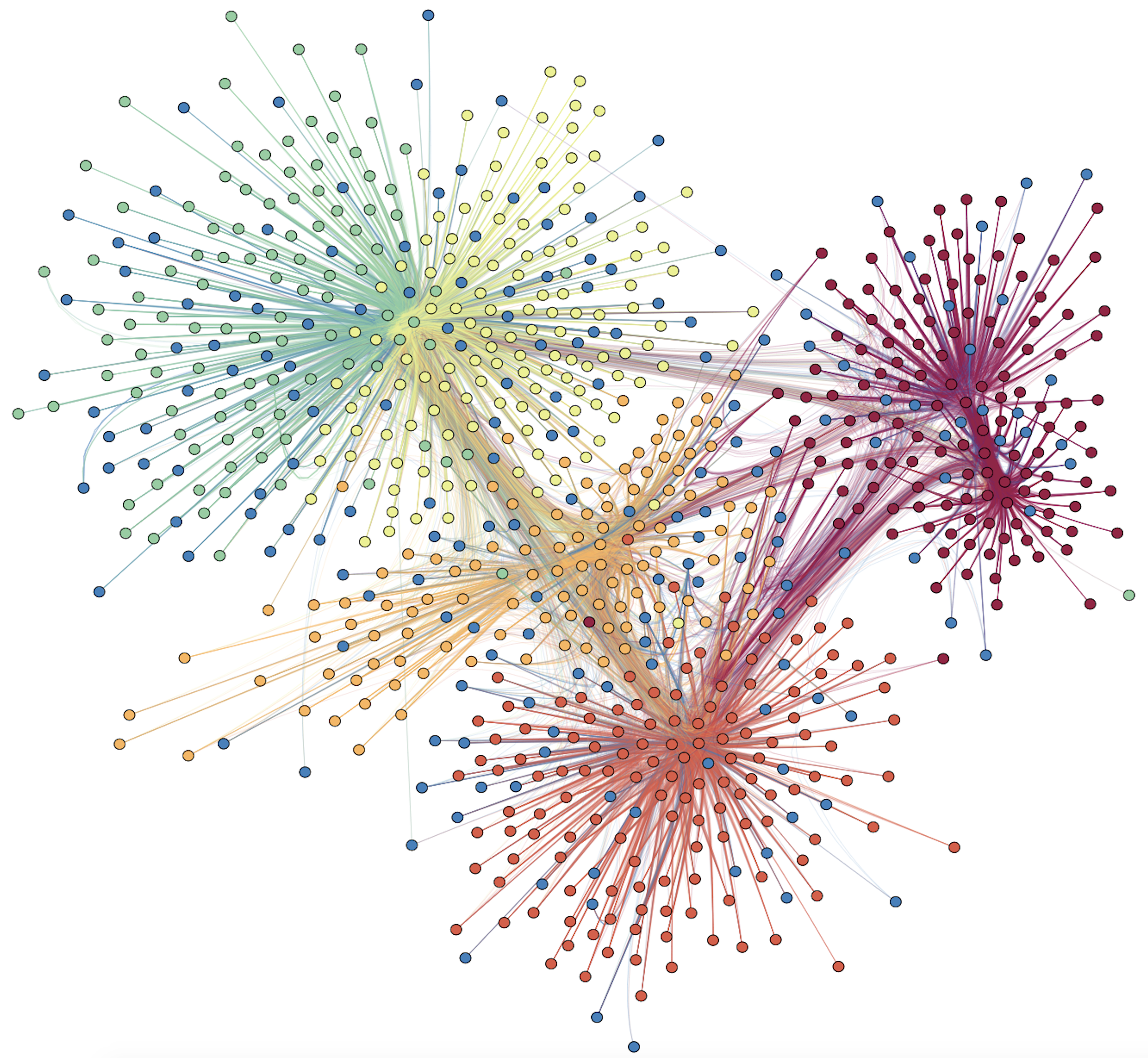}
 \end{center}
Each dot represents a \textit{class} of elliptic or parabolic pdes and each path is a compactness theorem.

For instance, the study of all Poisson equations of the form $-\Delta u = f(x)$, say with $f \in L^p$, is one single dot. Caccioppoli-type energy estimates yield a connecting path from such dot to the (sought-after) dot representing all harmonic functions. Namely, if $u_n$ is a sequence of functions, say bounded in $L^2$, satisfying $- \Delta u_n = f_n(x)$, if $f_n \to 0$ in $L^p$, then, up to a subsequence, $u_n \to h$ and $h$ is harmonic.

The abstract concept of Tangent postulates that two connected PDEs should share an underlying regularity theory. The \textit{caliber} of the path determines how much of the regularity one can bring from one model to another. In the above example, if 
$$
	\Delta u = f(x), \quad f \in L^p,
$$
then for any $0< \lambda \ll 1$, the function $u_\lambda(x) := \dfrac{1}{\lambda^{2-\frac{n}{p}}} u(\lambda x)$
verifies 
$$\Delta u_\lambda = f_\lambda(x),$$
where
 $$
 	f_\lambda(x) = \lambda^{\frac{n}{p}} f(\lambda x).
$$	
One easily verifies that 
$$
	\|f_\lambda(x)\|_p \le \|f\|_p.
$$	
This means that through the path joining the Poisson equation, $-\Delta u = f(x) \in L^p$, and  the Laplace equation, $\Delta h = 0$, one can transport estimates of order O$(r^{2-\frac{n}{p}})$. Such estimates ultimately yield $C^{0, 2-\frac{n}{p}}-$regularity, if $0 < 2-\frac{n}{p} < 1$, or $C^{1, 1-\frac{n}{p}}-$regularity if $1 < 2-\frac{n}{p} < 2$. As usual, the cases $p=n$ or $p=\infty$ are a bit tricker, as logarithm defects appear (see \cite{TT1}).

In recent years, methods from Geometric Tangential Analysis have been significantly enhanced, amplifying their range of application and providing a more user-friendly platform for advancing these endeavours (cf. \cite{TT1, TT2, TeiUrb14, TT3, PT, AraTeiUrb18, AraTeiUrb17}, to cite just a few). In what follows, we shall present a small sample of problems that can be tackled by methods coming from GTA.

\section {The parabolic $p-$Poisson equation}

As a first example, we consider the degenerate parabolic $p-$Poisson equation
\begin{equation} \label{eq_intro01}
	u_t - \mathrm{div} \left( |\nabla u|^{p-2} \nabla u \right) = f , \qquad p> 2,
\end{equation}

\noindent with a source term $f \in L^{q,r} (U_T) \equiv L^r(0,T;L^q(U))$, where the exponents satisfy the conditions
\begin{equation} \label{borderline1}
\frac{1}{r}+\frac{n}{pq} <1
\end{equation}
and
\begin{equation}\label{borderline2}
\frac{2}{r}+\frac{n}{q} >1.
\end{equation}
The first assumption is the standard minimal integrability condition that guarantees the existence of locally H\"older continuous weak solutions, while \eqref{borderline2} defines the borderline setting for optimal H\"older type estimates. For instance, when $r=\infty$,  conditions  \eqref{borderline1} and  \eqref{borderline2} enforce 
$$\dfrac{n}{p} < q < n,$$
which corresponds to the known range of integrability required in the elliptic theory for local $C^{0,\alpha}$ estimates to be available. 

We have shown in \cite{TeiUrb14} that weak solutions are locally of class $C^{0,\alpha}$ in space, with
$$\alpha := \frac{(pq-n)r-pq}{q[(p-1)r-(p-2)]} = \frac{p\left( 1-\displaystyle \frac{1}{r}-\frac{n}{pq}\right)}{\left( \displaystyle\frac{2}{r} + \frac{n}{q}-1\right) + p\left( 1-\displaystyle\frac{1}{r}-\frac{n}{pq}\right)},$$
a precise and sharp expression for the H\"older exponent in terms of $p$, the integrability of the source and the space dimension $n$. Observe that $0<\alpha<1$, in view of \eqref{borderline1} and \eqref{borderline2}.  

We also have that $u$ is of class $C^{0,\frac{\alpha}{\theta}}$ in time, where 
$$\theta := \alpha + p-(p-1)\alpha =  p - (p-2)\alpha = \alpha 2 + (1-\alpha) p$$
is the $\alpha-$interpolation between $2$ and $p$. If $p=2$, we have $\theta = 2$. For $p>2$, we have $2<\theta<p$, since $0 < \alpha < 1$.

The regularity proof develops along the following lines.

\bigskip

\begin{description}

\item[Step 1 - closeness to $p$-caloric] We first establish a key compactness result that states that if the source term $f$ has a small norm in $L^{q,r}$, then a solution $u$ to \eqref{eq_intro01} is close to a $p-$caloric function in an inner subdomain. The proof is by contradiction and uses compactness driven from a Caccioppoli-type energy estimate and a control of the time derivative due to Lindqvist in \cite{Lindq}.

\bigskip

\item[Step 2 - geometric iteration] Then, we explore the approximation by $p-$caloric functions and  the fact that $p-$caloric functions are \textit{universally} Lipschitz continuous in space and $C^{0, \frac{1}{2}}$ in time.  More precisely, we show there exist $\epsilon >0$ and $0< \lambda \ll 1/2$, depending only on $p$, $n$ and $\alpha$, such that if $\|f\|_{L^{q,r}(G_{1})} \leq \epsilon$ and $u$ is a local weak solution of \eqref{eq_intro01} in $G_{1}$, with $\|u\|_{p,avg,G_{1}}  \le 1$, then there exists a convergent sequence of real numbers $\{c_k\}_{k\ge 1}$, with
\begin{equation} \label{const}
|c_k - c_{k+1}| \leq c(n,p) \left( \lambda^\alpha \right)^k,
\end{equation}
such that
\begin{equation} \label{conclusion}
\|u-c_k\|_{p,avg,G_{\lambda^{k}}} \leq \left(\lambda^k\right)^\alpha,
\end{equation}
where the intrinsic $\theta-$parabolic cylinder is defined by
$$G_\tau := \left( - \tau^{\theta}, 0 \right) \times B_\tau(0), \quad \tau >0$$
and the averaged norm is 
$$\| v\|_{p, avg,Q} := \left( \, \intav{Q} |v|^p \, dx dt \right)^{1/p} = |Q|^{-1/p} \| v\|_{p,Q}.$$

\bigskip

\item[Step 3 - smallness regime] The smallness regime required is not restrictive since we can fall into that framework by scaling and contraction. Indeed, given a solution $u$, let 
$$v(x,t) = \varrho u ( x, \varrho^{p-2} t)$$
($\varrho$ to be fixed), which is a solution of \eqref{eq_intro01} with 
$$ \tilde{f} (x,t) =  \varrho^{p-1} f ( x, \varrho^{p-2} t).$$
Just choose $0<\varrho <1$ such that
$$\| v \|_{p,avg,G_{1}}^p \leq \varrho^{2} \|u\|_{p,avg,G_1}^p \leq 1$$
and
$$\| \tilde{f} \|_{L^{q,r} (G_1)}^r = \varrho^{(p-1) r - (p-2)} \| f \|_{L^{q,r} (G_1)}^r \leq \epsilon^r,$$
observing that, trivially, $(p-1) r - (p-2) >0$.

\bigskip

\item[Step 4 - H\"older via Campanato] Since the sequence $\{c_k\}$ is convergent, due to \eqref{const}, let
$$\bar{c} := \lim_{k \to \infty} c_k.$$
It follows from \eqref{conclusion} that, for arbitrary $0<r<\frac{1}{2}$,
$$ \intav{G_r} |u -\bar{c}  |^p \, dx dt  \leq C r^{p\alpha}.$$
Standard covering arguments and the characterisation of H\"older continuity of Campanato--Da Prato give the local $C^{0;\alpha, \alpha/\theta}$ -- continuity.

\end{description}

\bigskip

To highlight the extent to which our result is sharp, we project it into the state of the art of the theory. For the linear case $p=2$, we obtain
$$\alpha = 1 - \left( \frac{2}{r} + \frac{n}{q} -1 \right),$$
which is the optimal H\"older exponent for the non-homogeneous heat equation, and is in accordance with estimates obtained by energy considerations.
When $p\rightarrow \infty$, we have $\alpha \rightarrow 1^-$, which gives an indication of the expected locally Lipschitz regularity for the case of the parabolic infinity-Laplacian. When the source $f$ is independent of time, or else bounded in time, that is $r=\infty$, we obtain
$$\alpha = \frac{pq-n}{q(p-1)} = \frac{p}{p-1} \cdot \frac{q-\frac{n}{p}}{q},$$
which is exactly the optimal exponent of the elliptic case, easily obtained with the help of nonlinear Calder\'on-Zygmund theory and Morrey embeddings (see \cite{Ming}, and also \cite{T1} for borderline scenarios of integrability). 

Within the general theory of $p-$parabolic equations, our result reveals a surprising feature. From the applied point of view, it is relevant to know what is the effect on the diffusion properties of the model as we dim the exponent $p$. Na\"ive physical interpretations could indicate that the higher the value of $p$, the less efficient should the diffusion properties of the $p-$parabolic operator turn out to be, \textit{i.e.}, one should expect a less efficient smoothness effect of the operator. For instance, this is verified in the sharp regularity estimate for $p-$harmonic functions in the plane \cite{IM}. On the contrary, our estimate implies that for $p-$parabolic inhomogeneous equations, the H\"older regularity theory improves as $p$ increases. In fact, a direct computation shows
$$\text{sign}\left ( \partial_p \alpha(p,n,q,r) \right) = \text{sign} \left ( q(2-r) + nr \right ) =+1,$$
in view of standard assumptions on the integrability exponents of the source term.

\section {The porous medium equation}

There are two crucial differences with respect to the previous case when treating the porous medium equation (cf. \cite{Vaz07})
\begin{equation}
u_{t} - {\rm div} \left( m |u|^{m-1} \nabla u \right) = f, \quad m >1.
\label{pme}
\end{equation}
One is that adding a constant to a solution does not produce another solution, which somehow precludes the use of Campanato theory and requires a different technical approach to the H\"older regularity. The other is of a more fundamental nature, namely the hitherto unknown optimal regularity in the homogeneous case, leading to an extra dependence in the sharp H\"older exponent. Only for $n=1$,  it is proven in \cite{AroCaf86} that 
$$\alpha_0=\min \left\{1, \frac{1}{m-1}\right\} $$ 
but this is not the case in higher dimensions as corroborated by the celebrated counter-example in \cite{AroGra93}. 

A  locally bounded  function 
$$u \in C_{\rm loc} \left( 0,T ; L_{\rm loc}^{2}(U) \right), \qquad {\rm with} \quad |u|^{\frac{m+1}{2}} \in L_{\rm loc}^{2} \left( 0,T;W_{\rm loc}^{1,2}(U) \right)$$
is a local weak solution of \eqref{pme} if, for every compact set $K \subset U$ and every subinterval  $[t_{1}, t_{2}] \subset (0, T] $, we have
$$\left. \int_{K} u \varphi  \right|_{t_{1}}^{t_{2}} + \int_{t_{1}}^{t_{2}} \int_{K} \left\{ -u\varphi_{t} + m |u|^{m-1}\nabla u \cdot \nabla \varphi \right\} =\int_{t_{1}}^{t_{2}} \int_{K} f \varphi ,$$
for all  test functions  
$$\varphi \in W_{\rm loc}^{1,2} \left( 0,T;L^{2}(K) \right) \cap  L_{\rm loc}^{2} \left( 0,T;W_{0}^{1,2}(K) \right).$$
It is clear that all integrals in the above definition are convergent (cf. \cite[\S 3.5]{DiBeGiaVes11}), interpreting the gradient term as
$$|u|^{m-1}\nabla u : =  \frac{2}{m+1}\, {\rm sign} (u) \, |u|^{\frac{m-1}{2}} \nabla |u|^{\frac{m+1}{2}}.$$

For a source term $f \in L^{q,r} (U_T) \equiv L^r(0,T;L^q(U))$, with
$$\frac{1}{r}+\frac{n}{2q} <1,$$
it was shown in \cite{AMU} that locally bounded weak solutions of \eqref{pme} are locally of class $C^{0,\gamma}$ in space, with
\begin{equation}\gamma =\frac{\alpha}{m}  , \qquad \alpha =  \min\left\{\alpha_{0}^-, \frac{m[(2q - n)r -2q]}{q[mr - (m-1)]}\right\} ,
\label{alfa}
\end{equation}
where $0< \alpha_{0} \leq 1$ denotes the optimal H\"older exponent for solutions of \eqref{pme} with $f \equiv 0$. This regularity class is to be interpreted in the following sense: if 
$$\frac{m[(2q - n)r -2q]}{q[mr - (m-1)]} < \alpha_0$$
then solutions are in $C^{0,\gamma}$, with 
$$\gamma = \frac{(2q - n)r -2q}{q[mr - (m-1)]};$$
if, alternatively, 
$$\frac{m[(2q - n)r -2q]}{q[mr - (m-1)]} \geq \alpha_0,$$
then solutions are in $C^{0,\gamma}$, for any $0<\gamma <\frac{\alpha_{0}}{m} $.

Observe that
$$\frac{m[(2q - n)r -2q]}{q[mr - (m-1)]} = \frac{2 m \left( 1- \displaystyle \frac{1}{r} - \frac{n}{2q}\right)}{\displaystyle m \left( 1-\frac{1}{r}\right) +  \frac{1}{r} } >0$$
and so indeed $\gamma >0$. Note also that 
$$\frac{m[(2q - n)r -2q]}{q[mr - (m-1)]} >1$$ 
if 
$$\left( 1 + \frac{1}{m}\right) \frac{1}{r}+\frac{n}{q} < 1,$$
and, as $q,r \rightarrow \infty$,
$$\frac{m[(2q - n)r -2q]}{q[mr - (m-1)]} \longrightarrow 2,$$
which means that after a certain integrability threshold it is the optimal regularity exponent of the homogeneous case that prevails, with
$$\alpha = \alpha_0^- \qquad {\rm and}  \qquad \gamma <\frac{\alpha_{0}}{m}<1.$$ 

The $C^{0,\frac{\gamma}{\theta}}$ regularity in time is also obtained in \cite{AMU}, with 
$$\theta=2 - \left(1-\frac{1}{m}\right)\alpha=  \alpha \left( 1+\frac{1}{m} \right) + \left( 1-\alpha \right) 2 $$ 
being the $\alpha-$interpolation between $1+\frac{1}{m}$ and $2$. Observe that for $m=1$ we obtain
$$\gamma = 1 - \left( \frac{2}{r} + \frac{n}{q} -1 \right) \qquad {\rm and} \qquad \theta = 2,$$
recovering the optimal H\"older regularity for the non-homogeneous heat equation, in accordance with estimates obtained by energy considerations. 

It is worth stressing that, as in the case of the $p$-Laplace equation, the integrability in time (respectively, in space) of the source affects the regularity in space (respectively, in time) of the solution. Here is a snapshot of the regularity proof.

\bigskip

\begin{description}

\item[Step 1 - a Caccioppoli estimate]  An equivalent definition of weak solution involving the Steklov average is instrumental in obtaining the following Caccioppoli estimate (cf. \cite[\S 3.6]{DiBeGiaVes11}), for a constant $C$, depending only on $n, m$ and  $K \times [t_{1}, t_{2}] $:
$$\qquad \sup_{t_{1}< t < t_{2}} \int_{K}u^{2}\xi^{2} + \int_{t_{1}}^{t_{2}}\int_{K} |u|^{m-1}\vert\nabla u\vert^{2}\xi^{2} $$
$$\qquad \leq  C\int_{t_{1}}^{t_{2}}\int_{K} u^{2}\xi \left|\xi_{t}\right|+\int_{t_{1}}^{t_{2}}\int_{K} |u|^{m+1} \left( \vert\nabla \xi\vert^{2} +\xi^2 \right)+ C\Vert f \Vert^{2}_{L^{q ,r}},$$
for all $\xi \in C_{0}^{\infty}(K \times (t_{1}, t_{2}))$ such that $\xi \in [0, 1].$ 

\bigskip

\item[Step 2 - the intrinsic geometric setting] It is crucial we perform our analysis in the adequate geometric setting, reflecting the degeneracy in the pde in the scaling of the space-time cylinders where the oscillation is measured. Given $0 < \alpha \leq 1$, let
\begin{equation}
\label{theta}
\theta : = 2 - \left(1-\frac{1}{m}\right)\alpha
\end{equation}
and define the intrinsic $\theta $-parabolic cylinder as
$$G_{\rho} := \left( -\rho^\theta, 0 \right) \times B_{\rho}(0), \quad \rho > 0.$$
Note that when $m=1$, we obtain the standard parabolic cylinders reflecting the natural homogeneity between space and time for the heat equation. Observe also that $\theta$ satisfies the bounds
$$1+ \frac{1}{m} \leq \theta < 2,$$ 
which will be instrumental in the sequel.

\bigskip

\item[Step 3 - link with the homogeneous case] Another vital ingredient is a way to link solutions of the non-homogeneous problem with solutions of the homogeneous case. This type of statement follows from the available compactness for the problem and is proved by contradiction. Using the Caccioppoli estimate and the Arzel\`a-Ascoli theorem, we show that for every $\delta >0$, there exists  $0 < \epsilon \ll 1$ such that, if $\Vert f \Vert_{L^{q ,r}(G_{1})} \leq \epsilon$ and $u$ is a local weak solution of $(\ref{pme})$ in $G_{1}$, with $\Vert u \Vert_{\infty, G_{1}} \leq 1$, then there exists $\phi$ such that  
$$\phi_{t} - {\rm div} \left( m |\phi|^{m-1} \nabla \phi \right) = 0 \quad {\rm in}\ G_{1/2}$$
and
$$\Vert u -\phi \Vert_{\infty, G_{1/2}} \leq \delta .$$

\bigskip

\item[Step 4 - geometric iteration] This is the heart of the proof. We show that there exist  $\epsilon >0$ and $0 < \lambda \ll 1/2$, depending only on $m,n$ and $\alpha$, such that, if $\Vert f \Vert_{L^{q ,r}(G_{1})} \leq \epsilon$ and $u$ is a local weak solution of \eqref{pme} in $G_{1}$, with $\Vert u \Vert_{\infty, G_{1}} \leq 1,$  then, for every $k \in \N$,
$$\Vert u  \Vert_{\infty, G_{\lambda^{k}}} \leq  (\lambda^k)^{\gamma},$$
provided
$$|u(0,0)|\leq \frac{1}{4} \left( \lambda^k \right)^{\gamma}.$$
Recall that $\gamma$ was fixed in \eqref{alfa} and let us proceed by induction. We start with the case $k=1$. Take $0<\delta <1$, to be chosen later, and apply Step 3 to obtain $0 < \epsilon \ll 1$ and a solution $\phi$ of the homogeneous pme in $G_{1/2}$  such that
$$\Vert u -\phi \Vert_{\infty, G_{1/2}} \leq \delta.$$
Since $\phi$ is locally $C_{x}^{\alpha_{0}} \cap C_{t}^{\alpha_{0}/2} $, we obtain 
\begin{eqnarray*}
\underset{(x, t) \in G_{\lambda}}{\sup} \vert \phi(x, t) - \phi(0, 0) \vert \leq C\lambda^{\frac{\alpha_{0}}{m}},
\end{eqnarray*}
for $C > 1$ universal, where $\lambda \ll 1$ is still to be chosen. In fact, for $(x, t) \in G_{\lambda}$,
\begin{eqnarray*}
\vert \phi(x, t) - \phi(0, 0) \vert &\leq &  \vert \phi(x, t) - \phi(0, t) \vert + \vert \phi(0, t) - \phi(0, 0) \vert \\
&\leq & c_{1}\vert x - 0 \vert^{\alpha_{0}} +  c_{2}\vert t - 0 \vert^{\alpha_{0}/2}\\
&\leq & c_{1}\lambda^{\alpha_{0}} + c_{2}\lambda^{\frac{\theta}{2}\alpha_{0}}\\
&\leq & C\lambda^{\frac{\alpha_{0}}{m}},
\end{eqnarray*}
since $\theta \geq 1+ \frac{1}{m} >\frac{2}{m}$. We can therefore estimate 
\begin{eqnarray*}
\sup_{G_\lambda} |u| &\leq &  \sup_{G_{1/2}}  |u-\phi | + \sup_{G_\lambda} |\phi-\phi(0,0)| \nonumber \\
& & + |\phi(0,0)-u(0,0)| + |u(0,0)| \nonumber \\
&\leq & 2 \delta + C\lambda^{\frac{\alpha_{0}}{m}}+ \frac{1}{4} \lambda^{\gamma} \label{solskajer}
\end{eqnarray*}
and the result follows from the choices
$$\lambda = {\left(\frac{1}{4C}\right)^{\frac{m}{\alpha_{0} - \alpha}}} \qquad {\rm and} \qquad \delta= \frac{1}{4} \lambda^{\gamma} .$$
Now suppose the conclusion holds for $k$ and let's show it also holds for $k+1$. Due to $\eqref{theta}$, the function $v:G_{1} \rightarrow \R$  defined by
$$v(x,t) = \frac{u(\lambda^{k}x, \lambda^{k\theta}t)}{\lambda^{\gamma k}}$$
solves
$$v_{t}- {\rm div} \left( m |v|^{m-1}\nabla v \right)= \lambda^{k(2-\alpha)}f(\lambda^{k}x, \lambda^{k\theta}t)= \tilde{f}(x,t).$$
We have
\begin{eqnarray*}
\Vert \tilde{f} \Vert_{L^{q,r}(G_{1})}^{r} &=& \int_{-1}^{0} \left(\int_{B_{1}} \left\vert \tilde{f}(x,t) \right\vert^{q}  \right)^{r/q}\\
&=& \int_{-1}^{0}\left(  \int_{B_{1}} \lambda^{k(2 - \alpha )q}     \left\vert  f(\lambda^{k}x, \lambda^{k\theta}t) \right\vert^{q}  \right)^{r/q}\\
&=& \int_{-1}^{0}\left(  \int_{B_{\lambda^k}} \lambda^{k(2 - \alpha)q - kn} \left\vert  f(x, \lambda^{k\theta}t) \right\vert^{q}  \right)^{r/q}\\
&=&  \lambda^{[k(2 - \alpha )q - kn]\frac{r}{q} } \int_{-1}^{0}\left( \int_{B_{\lambda^{k}}}  \left\vert f(x, \lambda^{k\theta}t) \right\vert^{q}  \right)^{r/q}\\
&=&  \lambda^{[k(2 - \alpha )q - kn]\frac{r}{q} - k\theta}  \int_{-\lambda^{k\theta}}^{0}\left(  \int_{B_{\lambda^{k}}}  \left\vert f(x, t) \right\vert^{q}  \right)^{r/q}\\
\end{eqnarray*}
and, since
$$\left[ (2 - \alpha)q - n \right]\frac{r}{q} - \theta \geq  0$$
due to \eqref{alfa}, we get
\begin{eqnarray*}
\qquad \Vert \tilde{f} \Vert_{L^{q,r}(G_{1})} \leq \Vert f \Vert_{L^{q,r}( (-\lambda^{\theta k}, 0) \times B_{\lambda^{k}})} \leq \Vert f \Vert_{L^{q,r}(G_{1})} \leq \epsilon,
\end{eqnarray*}
which entitles $v$ to the case $k=1$. Note that $\Vert v \Vert_{\infty,G_{1}} \leq 1$, due to the induction hypothesis, and 
$$\left| v(0,0) \right| = \left| \frac{u(0,0)}{ \left( \lambda^k \right)^{\gamma}} \right| \leq \left| \frac{\frac{1}{4} \left( \lambda^{k+1} \right)^{\gamma}}{ \left( \lambda^k \right)^{\gamma}} \right| \leq \frac{1}{4} \lambda^{\gamma}.$$ 
It then follows that 
$$\Vert v \Vert_{\infty, G_{\lambda}} \leq \lambda^{\gamma},$$
which is the same as
$$\Vert u \Vert_{\infty, G_{\lambda^{k+1}}} \leq  \left( \lambda^{k+1} \right)^\gamma,$$
and the induction is complete. 

\bigskip

\item[Step 5 - the smallness regime] The final ingredient is showing that the smallness regime previously required is not restrictive. In fact, we can show that if $u$ is a local weak solution of \eqref{pme} in $G_{1}$ then, for every $0<r<\lambda$, 
\begin{equation}
\|u  \|_{\infty, G_{r}} \leq  C\, r^{\gamma},
\label{irala}
\end{equation}
provided
\begin{equation}
|u(0,0)|\leq \frac{1}{4} r^{\gamma}.
\label{santi}
\end{equation}
To see this, take
$$v(x,t) = \rho u \left( \rho^{a}x, \rho^{(m-1) + 2a}t \right) $$
with $\rho, a$ to be fixed. It solves
$$\qquad v_{t}- {\rm div} (m |v|^{m-1}\nabla v)= \rho^{m + 2a}f(\rho^{a}x, \rho^{(m-1)+ 2a)}t) = \tilde{f}(x,t)$$
and satisfies the bounds
$$\Vert v \Vert_{\infty,G_{1}} \leq  \rho \Vert u \Vert_{\infty,G_{1}}$$
and
$$\Vert \tilde{f} \Vert_{L^{q,r}(G_{1})}^{r} = \rho^{(m + 2a)r - a(n\frac{r}{q} + 2) - (m-1)}\Vert f \Vert_{L^{q,r}(G_{1})}^{r}.$$
Choosing $a > 0$ such that
$$(m + 2a)r - a \left( \frac{nr}{q} + 2 \right) - (m-1) > 0,$$
which is always possible, and $0 < \rho \ll 1$, we enter the smallness regime. Now, given $0<r<\lambda$, there exists $k \in \N$ such that 
$$\lambda^{k+1} < r \leq \lambda^{k}.$$
Since $|u(0,0)|\leq \frac{1}{4} r^{\gamma} \leq  \frac{1}{4} (\lambda^{k})^{\gamma}$, it follows that
$$\Vert u  \Vert_{\infty, G_{\lambda^{k}}} \leq  (\lambda^k)^{\gamma}$$
and, for $C=\lambda^{-\gamma}$,
$$\Vert u  \Vert_{\infty, G_{r}} \leq \Vert u  \Vert_{\infty, G_{\lambda^{k}}} \leq  (\lambda^k)^{\gamma} < \left( \frac{r}{\lambda} \right)^{\gamma} = C\, r^{\gamma} .$$

\bigskip

\item[Step 6 - sharp H\"older regularity] We finally show there exists a universal constant $K$ such that 
$$\|u-u(0,0)\|_{\infty, G_r} \leq K r^\gamma,$$
which is the $C^{0,\gamma}$ regularity at the origin.
We know, \textit{a priori}, that $u$ is continuous so 
$$\mu:= (4| u(0,0) |)^{1/\gamma} \geq 0$$
is well defined. Taking any radius $0<r<\lambda$, we have three alternative cases.

\medskip
\begin{itemize}

\item If $\mu \leq r<\lambda$, then \eqref{santi} holds and 
$$ \qquad \sup_{G_r} \left| u(x,t) - u(0,0)\right| \leq C\, r^{\gamma} + |u(0,0)| \leq \left( C+\frac{1}{4}\right) r^{\gamma}$$
follows from \eqref{irala}.
\medskip

\item If $0<r<\mu$, we consider the function 
$$w(x,t):=\frac{u(\mu x, \mu^\theta t)}{\mu^{\gamma}},$$
which solves a uniformly parabolic equation in $G_{\rho_0}$, for a radius $\rho_0$ depending only on the data. This gives an estimate that, written in terms of $u$, reads
$$\sup_{G_{r}} \left| u(x,t)  - u(0,0)  \right| \leq C\, r^{\gamma}, \quad \forall \, 0<r<\mu\frac{\rho_0}{2}.$$

\item Finally, for $\mu\frac{\rho_0}{2} \leq r < \mu$, we have
\begin{eqnarray*}
\qquad \sup_{G_{r}} \left| u(x,t)  - u(0,0)  \right| & \leq &   \sup_{G_{\mu}} \left| u(x,t)  - u(0,0)  \right| \nonumber\\
& \leq & C\, \mu^{\gamma} \leq C \left( \frac{2r}{\rho_0} \right)^{\gamma}= \tilde{C} r^{\gamma} .
\end{eqnarray*}

\end{itemize}

\end{description}

\section {The doubly nonlinear equation}

The methods and techniques of the previous two sections can be combined to obtain sharp regularity results for solutions of the inhomogeneous degenerate {\it doubly nonlinear parabolic equation}
$$u_t - \mathrm{div} \left( m |u|^{m-1}  |\nabla u|^{p-2} \nabla u \right) = f , \qquad p> 2 \quad m >1,$$
which appears, for example, in the contexts of non-Newtonian fluid dynamics, plasma physics, ground water problems or image processing. 

The local H\"older continuity of bounded weak solutions is established in \cite{Iva89, PorVes}. For a source term $f \in L^{q,r} (U_T) \equiv L^r(0,T;L^q(U))$, with 
$$\frac{1}{r}+\frac{n}{pq} <1\qquad {\rm and} \qquad \frac{3}{r}+\frac{n}{q} >2,$$
it is proven in \cite{Jani} that locally bounded weak solutions are locally of class $C^{0,\beta}$ in space with
$$\beta =\frac{\alpha (p-1)}{m+p-2}  , \qquad \alpha =  \min\left\{\alpha_{\ast}^-, \frac{(m+p-2)[(pq - n)r -pq]}{q(p-1)[(r-1)(m+p-2)+1]}\right\} ,$$
where $0< \alpha_{\ast} \leq 1$ denotes the optimal (unknown) H\"older exponent for solutions of the homogeneous case. The regularity class is to be interpreted as in the case of the porous medium equation.

Observe that when $m=1$, the equation becomes the degenerate parabolic $p-$Poisson equation, for which $\alpha_{\ast} = 1$, and we recover the exponent
$$\alpha := \frac{(pq-n)r-pq}{q[(p-1)r-(p-2)]}$$
of section 3. For $p=2$, we have the porous medium equation and obtain the exponent \eqref{alfa} of section 4.

\bigskip

{\small  \noindent{\bf Acknowledgments.} JMU partially supported by FCT -- Funda\c c\~ao para a Ci\^encia e a Tecnologia, I.P., through project PTDC/MAT-PUR/28686/2017, and by the Centre for Mathematics of the University of Coimbra -- UID/MAT/00324/2013, funded by the Portuguese government through FCT and co-funded by the European Regional Development Fund through Partnership Agreement PT2020.}

\end{document}